\documentclass[12pt]{amsart}

\usepackage{amsmath}
\usepackage{amssymb}
\usepackage{array}

\newtheorem{theorem}{Theorem}

\newtheorem{corollary}{Corollary}

\theoremstyle{definition}

\begin{document}

\title{Vacca-type  series for values of the generalized-Euler-constant function and its derivative}

\author{Kh.~Hessami Pilehrood$^{1}$}

\address{Institute for Studies in Theoretical Physics and Mathematics
(IPM), Tehran, Iran}

\curraddr{Mathemetics Department, Faculty of Basic Sciences,
Shahrekord University, Shahrekord, P.O. Box 115, Iran.}
\email{hessamik@ipm.ir, hessamik@gmail.com,
hessamit@ipm.ir}
\thanks{$^1$ This research was in part supported by a grant
from IPM (No. 86110025)}

\author{T.~Hessami Pilehrood$^2$}
\thanks{$^2$ This research was in part supported by a grant from
IPM (No. 86110020)}


\date{}


\begin{abstract}
We generalize well-known Catalan-type integrals for Euler's
constant to values of the ge\-ne\-ra\-li\-zed-Euler-constant
function and its derivatives.
 Using generating functions appeared in
these integral representations we give new Vacca and
Ramanujan-type series for values of the generalized-Euler-constant
function and Ad\-dison-type series for values of the
generalized-Euler-constant function and its derivative. As a
consequence, we get base $B$ rational series for
$\log\frac{4}{\pi},$ $\frac{G}{\pi}$ (where $G$ is Catalan's
constant), $\frac{\zeta'(2)}{\pi^2}$ and also for logarithms of
Somos's and Glaisher-Kinkelin's constants.

\end{abstract}

\maketitle

\section{Introduction}
\label{intro}

\vspace{0.3cm}

In \cite{s}, J.~Sondow proved the following two formulas:
\begin{equation}
\qquad\gamma=\sum_{n=1}^{\infty}\frac{N_{1,2}(n)+N_{0,2}(n)}{2n(2n+1)},
\label{eq1}
\end{equation}
\begin{equation}
\log\frac{4}{\pi}=\sum_{n=1}^{\infty}\frac{N_{1,2}(n)-N_{0,2}(n)}{2n(2n+1)},
\label{eq2}
\end{equation}
where $\gamma$ is Euler's constant and $N_{i,2}(n)$ is the number
of $i$'s in the binary expansion of $n.$ The series (\ref{eq1}) is
equivalent to  the well-known Vacca series \cite{va}
\begin{equation}
\gamma=\sum_{n=1}^{\infty}(-1)^n\frac{\lfloor\log_2n\rfloor}{n}=\sum_{n=1}^{\infty}
(-1)^n\frac{N_{1,2}\bigl(\lfloor\frac{n}{2}\rfloor\bigr)+
N_{0,2}\bigl(\lfloor\frac{n}{2}\rfloor\bigr)}{n} \label{eq3}
\end{equation}
and both series (\ref{eq1}) and (\ref{eq3}) may be derived from
Catalan's integral \cite{ca}
\begin{equation}
\gamma=\int_0^1\frac{1}{1+x}\sum_{n=1}^{\infty}x^{2^n-1}\,dx.
\label{eq4}
\end{equation}
To see this it  suffices to note that
$$
G(x)=\frac{1}{1-x}\sum_{n=0}^{\infty}x^{2^n}=\sum_{n=1}^{\infty}
(N_{1,2}(n)+N_{0,2}(n))x^n
$$
is a generating function of the sequence $N_{1,2}(n)+N_{0,2}(n),$
(see \cite[sequence A070939]{oeis}), which is the binary length of
$n,$ rewrite (\ref{eq4}) as
$$
\gamma=\int_0^1(1-x)\frac{G(x^2)}{x}\,dx
$$
and integrate the power series termwise. In view of the equality
\begin{equation*}
1=\int_0^1\sum_{n=1}^{\infty}x^{2^n-1}\,dx,
\end{equation*}
which is  easily verified by termwise integration, (\ref{eq4}) is
equivalent to the formula
\begin{equation}
\gamma=1-\int_0^1\frac{1}{1+x}\sum_{n=1}^{\infty}x^{2^n}\,dx
\label{eq05}
\end{equation}
obtained independently by Ramanujan (see \cite[Cor. 2.3]{bebo}).
Catalan's integral (\ref{eq05}) gives the following rational
series for $\gamma:$
\begin{equation}
\gamma=1-\int_0^1(1-x)G(x^2)\,dx=1-\sum_{n=1}^{\infty}
\frac{N_{1,2}(n)+N_{0,2}(n)}{(2n+1)(2n+2)}. \label{eq06}
\end{equation}
Averaging (\ref{eq1}), (\ref{eq06}) and (\ref{eq4}), (\ref{eq05}),
respectively, we get Addison's series for $\gamma$ \cite{ad}
$$
\gamma=\frac{1}{2}+\sum_{n=1}^{\infty}\frac{N_{1,2}(n)+N_{0,2}(n)}%
{2n(2n+1)(2n+2)}
$$
and its corresponding integral
\begin{equation}
\gamma=\frac{1}{2}+\frac{1}{2}\int_0^1\frac{1-x}{1+x}\sum_{n=1}^{\infty}
x^{2^n-1}, \label{eq065}
\end{equation}
 respectively. Integrals (\ref{eq05}),
(\ref{eq4}) were generalized to an arbitrary integer base $B>1$ by
S.~Ramanujan and B.~C.~Berndt and D.~C.~Bowman (see \cite{bebo})
\begin{equation}
\gamma
=1-\int_0^1\left(\frac{1}{1-x}-\frac{Bx^{B-1}}{1-x^B}\right)
\sum_{n=1}^{\infty}x^{B^n}\,dx \qquad\qquad (\text{Ramanujan}),
\label{eq07}
\end{equation}
 \begin{equation}
  \gamma
=\int_0^1\left(\frac{B}{1-x^B}-\frac{1}{1-x}\right)\sum_{n=1}^{\infty}
x^{B^n-1} \qquad\qquad (\text{Berndt-Bowman}). \label{eq08}
\end{equation}
Formula (\ref{eq08}) implies the generalized Vacca series for
$\gamma$ (see \cite[Th. 2.6]{bebo}) proposed by L.~Carlitz
\cite{carl}
\begin{equation}
\gamma=\sum_{n=1}^{\infty}\frac{\varepsilon(n)}{n}\lfloor\log_Bn\rfloor,
\label{eq09}
\end{equation}
 where
\begin{equation}
\varepsilon(n)=\begin{cases}
B-1     & \quad \text{if}\quad B\quad \text{divides}\quad n \\
 -1 & \quad \text{otherwise},
 \end{cases}
 \label{eq095}
\end{equation}
and the averaging integral of (\ref{eq07}) and (\ref{eq08})
produces the generalized Addison series for $\gamma$ found by
Sondow in  \cite{s}
\begin{equation}
\gamma=\frac{1}{2}+\sum_{n=1}^{\infty}\frac{\lfloor\log_BBn\rfloor
P_B(n)}{Bn(Bn+1)\cdots (Bn+B)}, \label{eq010}
\end{equation}
 where $P_B(x)$ is a
polynomial of degree $B-2$ denoted by
\begin{equation}
P_B(x)=(Bx+1)(Bx+2)\ldots
(Bx+B-1)\sum_{m=1}^{B-1}\frac{m(B-m)}{Bx+m}. \label{polyn}
\end{equation}
 In this
short note, we generalize Catalan-type integrals (\ref{eq07}),
(\ref{eq08}) to values of the ge\-ne\-ra\-li\-zed-Euler-constant
function
\begin{equation}
\gamma_{a,b}(z)=\sum_{n=0}^{\infty}\left(\frac{1}{an+b}-
\log\left(\frac{an+b+1}{an+b}\right)\right)z^n, \qquad a,b
\in{\mathbb N}, \label{eq7}
\end{equation}
and its derivatives,  which is related to constants (\ref{eq1}),
(\ref{eq2}) as $\gamma_{1,1}(1)=\gamma,$
$\gamma_{1,1}(-1)=\log\frac{4}{\pi}.$
 Using generating functions appeared in
these integral representations we give new Vacca and
Ramanujan-type series for values of $\gamma_{a,b}(z)$ and
Addison-type series for values of $\gamma_{a,b}(z)$ and its
derivative. As a consequence, we get base $B$ rational series for
$\log\frac{4}{\pi},$ $\frac{G}{\pi},$ (where $G$ is Catalan's
constant), $\frac{\zeta'(2)}{\pi^2}$ and also for logarithms of
Somos's and Glaisher-Kinkelin's constants. We also mention on
connection of our approach to summation of series of the form
$$\
\sum_{n=1}^{\infty}N_{\omega,B}(n)Q(n,B) \quad\text{and}\quad
\sum_{n=1}^{\infty}\frac{N_{\omega,B}(n)P_B(n)}{Bn(Bn+1)\cdots
(Bn+B)},
$$
where $Q(n,B)$ is a rational function of $B$ and $n$
\begin{equation}
Q(n,B)=\frac{1}{Bn(Bn+1)}+\frac{2}{Bn(Bn+2)}+\cdots+\frac{B-1}{Bn(Bn+B-1)},
\label{eq5}
\end{equation}
and $N_{\omega,B}(n)$ is the number of occurrences of a word
$\omega$ over the alphabet $\{0,1,\ldots, B-1\}$ in the $B$-ary
expansion of  $n,$ considered in \cite{assh}. In this notation,
the generalized Vacca series (\ref{eq09}) can be written as
follows:
\begin{equation}
\gamma= \sum_{k=1}^{\infty}L_B(k)Q(k,B), \label{eq6}
\end{equation}
where $L_B(k):=\lfloor\log_BBk\rfloor
=\sum_{\alpha=0}^{B-1}N_{\alpha,B}(k)$ is the $B$-ary length of
$k.$ Indeed,  representing $n=Bk+r,$ $0\le r\le B-1$ and summing
in (\ref{eq09}) over $k\ge 1$ and $0\le r\le B-1$ we get
$$
\gamma=\sum_{k=1}^{\infty}\lfloor\log_BBk\rfloor\left(
\frac{B-1}{Bk}-\frac{1}{Bk+1}-\cdots-\frac{1}{Bk+B-1}\right)
=\sum_{k=1}^{\infty}\lfloor\log_BBk\rfloor Q(k,B).
$$
 By the same notation, the generalized
Addison series (\ref{eq010}) gives another base $B$ expansion of
Euler's constant
\begin{equation}
\gamma=\frac{1}{2}+\sum_{n=1}^{\infty}
\frac{L_B(n)P_B(n)}{Bn(Bn+1)\cdots (Bn+B)}
=\frac{1}{2}+\sum_{n=1}^{\infty}L_B(n)
\left(Q(n,B)-\frac{B-1}{2Bn(n+1)}\right) \label{eu}
\end{equation}
 which converges faster than (\ref{eq6}) to $\gamma.$ Here we used
the fact that
$$
\sum_{n=1}^{\infty}\sum_{\alpha=0}^{B-1}\frac{N_{\alpha,
B}(n)}{n(n+1)}=\frac{B}{B-1},
$$
which can be easily checked by \cite[Section 3]{as}. On the other
hand,
\begin{equation*}
\begin{split}
Q&(n,B)-\frac{B-1}{2Bn(n+1)}=\frac{1}{2}\sum_{m=1}^{B-1}\left(\frac{1}{Bn}-
\frac{2}{Bn+m}+\frac{1}{Bn+B}\right) \\
&=\frac{1}{Bn(Bn+B)}\sum_{m=1}^{B-1}
\left(2m-B+\frac{2m(B-m)}{Bn+m}\right)=\frac{P_B(n)}{Bn(Bn+1)\cdots(Bn+B)}.
\end{split}
\end{equation*}
{\bf \small Acknowledgements:} Both authors thank the Max Planck
Institute for Mathematics at Bonn where this research was carried
out. Special gratitude is due to professor B.~C.~Berndt for
providing paper \cite{bebo}.


\section{Analytic continuation}

\vspace{0.3cm}

We consider the generalized-Euler-constant function
$\gamma_{a,b}(z)$ defined in (\ref{eq7}), where $a, b$ are
positive real numbers, $z\in {\mathbb C},$ and the series
converges when $|z|\le 1.$ We show that $\gamma_{a,b}(z)$ admits
an analytic continuation to the domain ${\mathbb C}\setminus
[1,+\infty).$ The following theorem is a slight modification of
\cite[Th.3]{sh}.
\begin{theorem}
Let $a, b$ be positive real numbers, $z\in {\mathbb C},$ $|z|\le
1.$ Then
\begin{equation}
\gamma_{a,b}(z)=\int_0^1\int_0^1\frac{(xy)^{b-1}(1-x)}{(1-zx^ay^a)(-\log
xy)}\,dxdy=\int_0^1\frac{x^{b-1}(1-x)}{1-zx^a}\left(\frac{1}{1-x}+
\frac{1}{\log x}\right)\,dx. \label{eq8}
\end{equation}
The integrals converge for all $z\in {\mathbb C}\setminus
(1,+\infty)$ and give the analytic continuation of the
generalized-Euler-constant function $\gamma_{a,b}(z)$ for $z\in
{\mathbb C}\setminus [1,+\infty).$
\end{theorem}
{\bf Proof.} Denoting the double integral in (\ref{eq8}) by $I(z)$
and for $|z|\le 1,$ expanding $(1-zx^ay^a)^{-1}$ in a geometric
series we have
\begin{equation*}
\begin{split}
I(z)&=\sum_{k=0}^{\infty}z^k\int_0^1\int_0^1\frac{(xy)^{ak+b-1}(1-x)}{(-\log
xy)}\,dxdy \\
&=\sum_{k=0}^{\infty}z^k\int_0^1\int_0^1\int_0^{+\infty}(xy)^{t+ak+b-1}(1-x)\,dxdydt
\\
&=\sum_{k=0}^{\infty}z^k\int_0^{+\infty}\left(\frac{1}{(t+ak+b)^2}-\Bigl(
\frac{1}{t+ak+b}-\frac{1}{t+ak+b+1}\Bigr)\right)\,dt=\gamma_{a,b}(z).
\end{split}
\end{equation*}
On the other hand, making the change of variables $u=x^a,$ $v=y^a$
in the double integral we get
$$
I(z)=\frac{1}{a}\int_{0}^1\int_{0}^1\frac{(uv)^{\frac{b}{a}-1}(1-u^{\frac{1}{a}})}%
{(1-zuv)(-\log uv)}\,dudv.
$$
Now by \cite[Corollary 3.3]{gs}, for $z\in {\mathbb C}\setminus
[1,+\infty)$ we have
$$
I(z)=\frac{1}{a}\Phi\Bigl(z,1,\frac{b}{a}\Bigr)-\frac{\partial\Phi}{\partial
s}\Bigl(z,0,\frac{b}{a}\Bigr)+\frac{\partial\Phi}{\partial
s}\Bigl(z,0,\frac{b+1}{a}\Bigr),
$$
where $\Phi(z,s,u)$ is the Lerch transcendent, a holomorphic
function in $z$ and $s,$ for $z\in {\mathbb C}\setminus
[1,+\infty)$ and all complex $s$ (see \cite[Lemma 2.2]{gs}), which
is the analytic continuation of the series
$$
\Phi(z,s,u)=\sum_{n=0}^{\infty}\frac{z^n}{(n+u)^s}, \qquad u>0.
$$
To prove the second equality in (\ref{eq8}), make the change of
variables $X=xy,$ $Y=y$ and integrate with respect to $Y.$ \qed

\begin{corollary} \label{c1}
Let $a,b$ be positive real numbers, $l\in {\mathbb N},$ $z\in
{\mathbb C}\setminus [1,+\infty).$ Then for the $l$-th derivative
we have
$$
\gamma_{a,b}^{(l)}(z)=\int_0^1\int_0^1\frac{(xy)^{al+b-1}(x-1)}%
{(1-zx^ay^a)^{l+1}\log xy}\,dxdy=
\int_0^1\frac{x^{la+b-1}(1-x)}{(1-zx^a)^{l+1}}
\left(\frac{1}{1-x}+\frac{1}{\log x}\right)\,dx.
$$
\end{corollary}

From Corollary \ref{c1}, \cite[Cor.3.3, 3.8, 3.9]{gs} and
\cite[Lemma 4]{assh} we get

\begin{corollary} \label{c1.1}
Let $a, b$ be positive real numbers, $z\in {\mathbb C}\setminus
[1, +\infty).$ Then the following equalities are valid:

$$
\gamma_{a,b}(1)=\log\Gamma\Bigl(\frac{b+1}{a}\Bigr)-\log\Gamma\Bigl(
\frac{b}{a}\Bigr)-\frac{1}{a}\psi\Bigl(\frac{b}{a}\Bigr),
$$

$$
\gamma_{a,b}(z)=\frac{1}{a}\Phi\Bigl(z,1,\frac{b}{a}\Bigr)-
\frac{\partial\Phi}{\partial s}\Bigl(z,0,\frac{b}{a}\Bigr)
+\frac{\partial\Phi}{\partial s}\Bigl(z,0,\frac{b+1}{a}\Bigr),
$$

\begin{equation*}
\begin{split}
\gamma'_{a,b}(z)=&-\frac{b}{a^2}\Phi\Bigl(z,1,\frac{b}{a}+1\Bigr)+
\frac{1}{a(1-z)}+\frac{b}{a}\frac{\partial\Phi}{\partial
s}\Bigl(z,0,\frac{b}{a}+1\Bigr)-  \frac{\partial\Phi}{\partial
s}\Bigl(z,-1,\frac{b}{a}+1\Bigr)-
\\[10pt]
&\frac{b+1}{a}\frac{\partial\Phi}{\partial
s}\Bigl(z,0,\frac{b+1}{a}+1\Bigr)+\frac{\partial\Phi}{\partial s}
\Bigl(z,-1,\frac{b+1}{a}+1\Bigr),
\end{split}
\end{equation*}

\noindent where $\Phi(z,s,u)$ is the Lerch transcendent and
$\psi(x)=\frac{d}{dx}\log\Gamma(x)$ is the logarithmic derivative
of the gamma function.
\end{corollary}

\vspace{0.5cm}

\section{Catalan-type integrals for $\gamma_{a,b}^{(l)}(z).$}

\vspace{0.3cm}

In \cite{bebo} it was demonstrated that for $x>0$ and any integer
$B>1,$ one has
\begin{equation}
\frac{1}{1-x}+\frac{1}{\log x}=\sum_{k=1}^{\infty}\frac{(B-1)+
(B-2)x^{\frac{1}{B^k}}+(B-3)x^{\frac{2}{B^k}}+\cdots+x^{\frac{B-2}{B^k}}}%
{B^k(1+x^{\frac{1}{B^k}}+x^{\frac{2}{B^k}}+\cdots+x^{\frac{B-1}{B^k}})}.
\label{eq9}
\end{equation}
The special cases $B=2,3$ of this equality can be found in
Ramanujan's third note book \cite[p.364]{ra}. Using this key
formula we prove the following generalization of integral
(\ref{eq08}).
\begin{theorem} \label{t2}
 Let $a, b, B>1$ be positive integers, $l$ a non-negative integer.
 If either $z\in {\mathbb C}\setminus [1,+\infty)$ and $l\ge 1,$
 or $z\in {\mathbb C}\setminus (1,+\infty)$ and $l=0,$ then
 \begin{equation}
 \gamma_{a,b}^{(l)}(z)=\int_0^1\left(\frac{B}{1-x^B}-\frac{1}{1-x}\right)
 F_l(z,x)\,dx
 \label{eq10}
 \end{equation}
 where
 \begin{equation}
 F_l(z,x)=\sum_{k=1}^{\infty}\frac{x^{(b+al)B^k-1}(1-x^{B^k})}%
 {(1-zx^{aB^k})^{l+1}}.
 \label{eq017}
 \end{equation}
 \end{theorem}

{\bf Proof.}  First we note that the series of variable $x$ on the
right-hand side of (\ref{eq9}) uniformly converges on $[0,1],$
since the absolute value of its general term does not exceed
$\frac{B-1}{2B^{k-1}}.$ Then for $l\ge 0,$  multiplying both sides
of (\ref{eq9}) by $\frac{x^{la+b-1}(1-x)}{(1-zx^a)^{l+1}}$ and
integrating over $0\le x\le 1$ we get
$$
\gamma_{a,b}^{(l)}(z)=\sum_{k=1}^{\infty}\int_0^1\frac{x^{la+b-1}(1-x)}{(1-zx^a)^{l+1}}
\cdot\frac{(B-1)+(B-2)x^{\frac{1}{B^k}}+\cdots+x^{\frac{B-2}{B^k}}}%
{B^k(1+x^{\frac{1}{B^k}}+x^{\frac{2}{B^k}}+\cdots+x^{\frac{B-1}{B^k}})}\,dx.
$$
Replacing $x$ by $x^{B^k}$ in each integral we find
\begin{equation*}
\begin{split}
\gamma_{a,b}^{(l)}(z)&=\sum_{k=1}^{\infty}\int_0^1\frac{x^{(la+b)B^k-1}(1-x^{B^k})}%
{(1-zx^{aB^k})^{l+1}}\cdot\frac{(B-1)+(B-2)x+\cdots+x^{B-2}}{1+x+x^2+\cdots+x^{B-1}}\,dx
\\ &=\int_0^1\left(\frac{B}{1-x^B}-\frac{1}{1-x}\right)F_l(z,x)\,dx,
\end{split}
\end{equation*}
as required. \qed

From Theorem \ref{t2} we readily get a generalization of
Ramanujan's integral.
\begin{corollary} \label{c2}
Let $a, b, B>1$ be positive integers, $l$ a non-negative integer.
 If either $z\in {\mathbb C}\setminus [1,+\infty)$ and $l\ge 1,$
 or $z\in {\mathbb C}\setminus (1,+\infty)$ and $l=0,$ then
\begin{equation}
\gamma_{a,b}^{(l)}(z)=\int_0^1\frac{x^{b+al-1}(1-x)}{(1-zx^a)^{l+1}}\,dx
+\int_0^1\left(\frac{Bx^{B}}{1-x^B}-\frac{x}{1-x}\right)F_l(z,x)\,dx.
\label{eq018}
\end{equation}
\end{corollary}

 {\bf Proof.} First we note that the series (\ref{eq017}) considered as a sum of
 functions of variable $x$ uniformly converges on $[0,
 1-\varepsilon]$ for any $\varepsilon>0.$  Then integrating
 termwise we have
 $$
 \int_0^{1-\varepsilon}F_l(z,x)\,dx=\sum_{k=1}^{\infty}
 \int_0^{1-\varepsilon}\frac{x^{(b+al)B^k-1}(1-x^{B^k})}%
 {(1-zx^{aB^k})^{l+1}}\,dx.
 $$
Making the change of variable $y=x^{B^k}$ in each integral we get
$$
\int_0^{1-\varepsilon}F_l(z,x)\,dx=\sum_{k=1}^{\infty}
\frac{1}{B^k} \int_0^{(1-\varepsilon)^{B^k}}\frac{y^{b+al-1}(1-y)}%
 {(1-zy^a)^{l+1}}\,dy.
 $$
Since the last series of variable $\varepsilon$ uniformly
converges on $[0,1],$ letting $\varepsilon$ tend to zero we get
\begin{equation}
\int_0^1F_l(z,x)\,dx=\frac{1}{B-1}\int_0^1
\frac{y^{b+al-1}(1-y)}{(1-zy^a)^{l+1}}\,dy. \label{eq019}
\end{equation}
Now from (\ref{eq10}) and (\ref{eq019}) it follows that
$$
\gamma_{a,b}^{(l)}(z)-\int_0^1
\frac{y^{b+al-1}(1-y)}{(1-zy^a)^{l+1}}\,dy= \int_0^1\left(
\frac{Bx^B}{1-x^B}-\frac{x}{1-x}\right)F_l(z,x)\,dx,
$$
and the proof is complete. \qed

Averaging both formulas (\ref{eq10}),  (\ref{eq018}) we get the
following generalization of integral (\ref{eq065}).
\begin{corollary} \label{c3}
Let $a, b, B>1$ be positive integers, $l$ a non-negative integer.
 If either $z\in {\mathbb C}\setminus [1,+\infty)$ and $l\ge 1,$
 or $z\in {\mathbb C}\setminus (1,+\infty)$ and $l=0,$ then
$$
\gamma_{a,b}^{(l)}(z)=\frac{1}{2}\int_0^1\frac{x^{b+al-1}(1-x)}{(1-zx^a)^{l+1}}\,dx
+\frac{1}{2}\int_0^1\left(\frac{B(1+x^B)}{1-x^B}-\frac{1+x}{1-x}\right)F_l(z,x)\,dx.
$$
\end{corollary}

\vspace{0.5cm}

\section{Vacca-type series for $\gamma_{a,b}(z)$ and $\gamma'_{a,b}(z).$}

\vspace{0.3cm}

\begin{theorem} \label{t3}
Let $a,b,B>1$ be positive integers,  $z\in {\mathbb C},$ $|z|\le
1.$   Then for
 the generalized-Euler-constant function $\gamma_{a,b}(z)$, the
 following expansion is valid:
 $$
 \gamma_{a,b}(z)=\sum_{k=1}^{\infty}a_kQ(k,B)=
 \sum_{k=1}^{\infty}a_{\lfloor\frac{k}{B}\rfloor}\frac{\varepsilon(k)}{k},
 $$
 where $Q(k,B)$ is a rational function given by {\rm (\ref{eq5})},
  $\{a_k\}_{k=0}^{\infty}$ is a sequence defined by the
 generating function
 \begin{equation}
 G(z,x)=\frac{1}{1-x}\sum_{k=0}^{\infty}\frac{x^{bB^k}(1-x^{B^k})}%
 {1-zx^{aB^k}}=\sum_{k=0}^{\infty}a_kx^k
 \label{eq11}
 \end{equation}
 and $\varepsilon(k)$ is denoted in {\rm (\ref{eq095})}.
 \end{theorem}
 {\bf Proof.} For $l=0,$ rewrite (\ref{eq10}) in the form
$$
\gamma_{a,b}(z)=\int_0^1\frac{1-x^B}{x}\left(\frac{B}{1-x^B}-\frac{1}{1-x}\right)
G(z,x^B)\,dx
$$
where $G(z,x)$ is defined in (\ref{eq11}). Then, since $a_0=0,$
we have
\begin{equation}
\gamma_{a,b}(z)=\int_0^1(B-1-x-x^2-\cdots-x^{B-1})\sum_{k=1}^{\infty}a_kx^{Bk-1}\,
dx. \label{eq12}
\end{equation}
Expanding $G(z,x)$ in a power series of $x$
$$
G(z,x)=\sum_{k=0}^{\infty}\sum_{m=0}^{\infty}z^mx^{(am+b)B^k}(1+x+\cdots+x^{B^k-1})
$$
we see that $a_k=O(\ln_B k).$ Therefore, by termwise integration
in (\ref{eq12}), which can be easily justified  by the same way as
in the proof of Corollary \ref{c2}, we get
\begin{equation*}
\begin{split}
\gamma_{a,b}(z)
&=\sum_{k=1}^{\infty}a_k\int_0^1[(x^{Bk-1}-x^{Bk})+(x^{Bk-1}-x^{Bk+1})+\cdots+
(x^{Bk-1}-x^{Bk+B-2})]\,dx \\
&=\sum_{k=1}^{\infty}a_kQ(k,B). \qquad \qed
\end{split}
\end{equation*}

\begin{theorem}
Let $a,b,B>1$ be positive integers,  $z\in {\mathbb C},$ $|z|\le
1.$   Then for
 the generalized-Euler-constant function, the
 following expansion is valid:
 $$
 \gamma_{a,b}(z)=\int_0^1\frac{x^{b-1}(1-x)}{1-zx^a}\,dx-\sum_{k=1}^{\infty}a_k
 \widetilde{Q}(k,B),
 $$
 where
 \begin{equation*}
 \begin{split}
 \widetilde{Q}(k,B)&=\frac{B-1}{Bk(k+1)}-Q(k,B) \\
&= \frac{B-1}{(Bk+B)(Bk+1)}+\frac{B-2}{(Bk+B)(Bk+2)}+\cdots+
 \frac{1}{(Bk+B)(Bk+B-1)}
 \end{split}
 \end{equation*}
  and the sequence
 $\{a_k\}_{k=1}^{\infty}$ is  defined in Theorem {\rm\ref{t3}}.
\end{theorem}

{\bf Proof.} From Corollary \ref{c2} with $l=0,$ by the same way
as in the proof of Theorem \ref{t3},
 we get
\begin{equation*}
\begin{split}
&\int_0^1\left(\frac{Bx^B}{1-x^B}-\frac{x}{1-x}\right)F_0(z,x)=
\int_0^1\frac{1-x^B}{x}\left(\frac{Bx^B}{1-x^B}-\frac{x}{1-x}\right)G(z,x^B)\,dx
\\
&=\int_0^1(Bx^{B-1}-(1+x+\cdots+x^{B-1}))\sum_{k=1}^{\infty}a_kx^{Bk}\,dx
\\ &=\sum_{k=1}^{\infty}a_k\int_0^1[(x^{Bk+B-1}-x^{Bk+B-2})+\cdots+
(x^{Bk+B-1}-x^{Bk+1})+(x^{Bk+B-1}-x^{Bk})]\,dx \\
&=-\sum_{k=1}^{\infty}a_k\widetilde{Q}(k,B). \qquad \qed
\end{split}
\end{equation*}

\begin{theorem} \label{t5}
Let $a,b,B>1$ be positive integers,  $z\in {\mathbb C},$ $|z|\le
1.$   Then for
 the generalized-Euler-constant function $\gamma_{a,b}(z)$ and its derivative, the
 following expansion is valid:
$$
\gamma_{a,b}^{(l)}(z)=\frac{1}{2}\int_0^1\frac{x^{b+al-1}(1-x)}{(1-zx^a)^{l+1}}\,dx
+\sum_{k=1}^{\infty}a_{k,l}\frac{P_B(k)}{Bk(Bk+1)\cdots(Bk+B)},
\qquad l=0,1,
$$
where $P_B(k)$ is a polynomial of degree $B-2$ given by {\rm
(\ref{polyn})}, $(z-1)^2+(l-1)^2\ne 0$ and   the sequence
$\{a_{k,l}\}_{k=0}^{\infty}$ is defined by the generating function
\begin{equation}
G_l(z,x)=\frac{1}{1-x}\sum_{k=0}^{\infty}\frac{x^{(b+al)B^k}(1-x^{B^k})}%
{(1-zx^{aB^k})^{l+1}}=\sum_{k=0}^{\infty}a_{k,l}x^k, \qquad l=0,
1. \label{eq24}
\end{equation}
\end{theorem}

{\bf Proof.} Expanding $G_l(z,x)$ in a power series of $x$
$$
G_l(z,x)=\sum_{k=0}^{\infty}\sum_{m=0}^{\infty}\binom{m+l}{l}z^mx^{(b+al+am)B^k}(1+x+x^2+
\cdots+x^{B^k-1})
$$
we see that $a_{k,l}=O(k^l\ln_Bk).$ Therefore, for $l=0,1,$ by
termwise integration we get
\begin{equation*}
\begin{split}
&\int_0^1\left(\frac{B(1+x^B)}{1-x^B}-\frac{1+x}{1-x}\right)F_l(z,x)dx
=\int_0^1\frac{1-x^B}{x}\left(\frac{B(1+x^B)}{1-x^B}-\frac{1+x}{1-x}\right)G_l(z,x^B)dx
\\
&=\int_0^1[(B-1)-2x-2x^2-\cdots-2x^{B-1}+(B-1)x^B]\sum_{k=1}^{\infty}a_{k,l}x^{Bk-1}\,dx
\\
&=\sum_{k=1}^{\infty}a_{k,l}\left(\frac{B-1}{Bk}-\frac{2}{Bk+1}-\frac{2}{Bk+2}-\cdots
-\frac{2}{Bk+B-1}+\frac{B-1}{Bk+B}\right) \\
&=2\sum_{k=1}^{\infty}a_{k,l}\frac{P_B(k)}{Bk(Bk+1)\cdots(Bk+B)},
\end{split}
\end{equation*}
where $P_B(k)$ is defined in (\ref{polyn}) and the last series
converges  since $\frac{P_B(k)}{Bk(Bk+1)\cdots(Bk+B)}=O(k^{-3}).$
Now our theorem easily follows from Corollary \ref{c3}. \qed

\vspace{0.5cm}

\section{Summation of series in terms of the Lerch transcendent}

\vspace{0.3cm}

It is easily seen that the generating function (\ref{eq24})
satisfies the following  functional equation:
\begin{equation}
G_l(z,x)-\frac{1-x^B}{1-x}G_l(z,x^B)=\frac{x^{b+al}}{(1-zx^a)^{l+1}},
 \label{eq25}
\end{equation}
which is equivalent to the  identity for series:
$$
\sum_{k=0}^{\infty}a_{k,l}x^k-(1+x+\cdots+x^{B-1})\sum_{k=0}^{\infty}a_{k,l}x^{Bk}
=\sum_{k=l}^{\infty}\binom{k}{l}z^{k-l}x^{ak+b}.
$$
Comparing coefficients of powers of $x$ we get an alternative
definition of the sequence $\{a_{k,l}\}_{k=0}^{\infty}$ by means
of the recursion
$$
a_{0,l}=a_{1,l}=\ldots=a_{al+b-1,l}=0
$$
and for $k\ge al+b,$
\begin{equation}
a_{k,l}=\begin{cases} a_{\lfloor\frac{k}{B}\rfloor,l} &
\qquad\text{if} \qquad
k\not\equiv b \pmod a, \\
a_{\lfloor\frac{k}{B}\rfloor,l}+\binom{(k-b)/a}{l}z^{\frac{k-b}{a}-l}
& \qquad\text{if} \qquad  k\equiv b\pmod a.
\end{cases}
\label{re}
\end{equation}
On the other hand, in view of Corollary \ref{c1.1},
$\gamma_{a,b}(z)$ and $\gamma'_{a,b}(z)$ can be explicitly
expressed in terms of the Lerch transcendent, $\psi$-function and
logarithm of the gamma function.  This allows us to sum the series
figured in Theorems \ref{t3}-\ref{t5} in terms of these functions.

\vspace{0.2cm}

\section{Examples of rational series}

\vspace{0.3cm}

{\bf Example 1.} Suppose that $\omega$ is a non-empty word over
the alphabet $\{0,1,\ldots,B-1\}.$ Then obviously $\omega$ is
uniquely defined by its length $|\omega|$ and its size
$v_B(\omega)$ which is the value of $\omega$ when interpreted as
an integer in base $B.$ Let $N_{\omega,B}(k)$ be the number of
(possibly overlapping) occurrences of the block $\omega$ in the
$B$-ary expansion of $k.$ Note that for every $B$ and $\omega,$
$N_{\omega,B}(0)=0,$ since the $B$-ary expansion of zero is the
empty word. If the word $\omega$ begins with $0,$ but
$v_B(\omega)\ne 0,$ then in computing $N_{\omega,B}(k)$ we assume
that the $B$-ary expansion of $k$ starts with an arbitrary long
prefix of $0$'s. If $v_B(\omega)=0$ we take for $k$ the usual
shortest $B$-ary expansion of $k.$

Now we consider equation (\ref{eq25}) with  $l=0,$ $z=1$
\begin{equation}
G(1,x)-\frac{1-x^B}{1-x}G(1,x^B)=\frac{x^b}{1-x^a} \label{eq26}
\end{equation}
  and for a given non-empty word $\omega,$ set in
(\ref{eq26}) $a=B^{|\omega|}$ and
\begin{equation*}
b=\begin{cases}
B^{|\omega|}     & \quad \text{if}\quad v_B(\omega)=0 \\
v_B(\omega) & \quad \text{if}\quad v_B(\omega)\ne 0.
 \end{cases}
\end{equation*}
Then by (\ref{re}), it is easily seen that
$a_k:=a_{k,0}=N_{\omega, B}(k),$ $k=1,2,\ldots,$ and by Theorem
\ref{t3}, we get one more proof of the following statement (see
\cite[Sections 3, 4.2]{assh}).
\begin{corollary} \label{c4}
Let $\omega$ be a non-empty word over the alphabet $\{0,1,\ldots,
B-1\}.$ Then
\begin{equation*}
\sum_{k=1}^{\infty}N_{\omega, B}(k)Q(k,B)=\begin{cases}
\gamma_{B^{|\omega|},v_B(\omega)}(1)     & \quad \text{if}\quad v_B(\omega)\ne 0 \\
\gamma_{B^{|\omega|}, B^{|\omega}|}(1) & \quad \text{if}\quad
v_B(\omega)= 0.
 \end{cases}
\end{equation*}
\end{corollary}
By Corollary \ref{c1.1}, the right-hand side of the last equality
can be calculated explicitly and we have
\begin{equation}
\sum_{k=1}^{\infty}N_{\omega, B}(k)Q(k,B)=\begin{cases}
\log\Gamma\left(\frac{v_B(\omega)+1}{B^{|\omega|}}\right)-
\log\Gamma\left(\frac{v_B(\omega)}{B^{|\omega|}}\right)-
\frac{1}{B^{|\omega|}}\psi\left(\frac{v_B(\omega)}{B^{|\omega|}}\right)     &  \text{if}\,\, v_B(\omega)\ne 0 \\
\log\Gamma\left(\frac{1}{B^{|\omega|}}\right)+
\frac{\gamma}{B^{|\omega|}}-|\omega|\log B &  \text{if}\,\,
v_B(\omega)= 0.
 \end{cases}
 \label{eq28}
\end{equation}
\begin{corollary} \label{c54}
Let $\omega$ be a non-empty word over the alphabet $\{0,1,\ldots,
B-1\}.$ Then
\begin{equation*}
\begin{split}
&\sum_{k=1}^{\infty}\frac{N_{\omega, B}(k)P_B(k)}{Bk(Bk+1)\cdots
(Bk+B)} \\ &=\begin{cases}
\gamma_{B^{|\omega|},v_B(\omega)}(1)-\frac{1}{2B^{|\omega|}}
\left(\psi\Bigl(\frac{v_B(\omega)+1}{B^{|\omega|}}\Bigr)
-\psi\Bigl(\frac{v_B(\omega)}{B^{|\omega|}}\Bigr)\right)     & \quad \text{if}\quad v_B(\omega)\ne 0 \\
\gamma_{B^{|\omega|}, B^{|\omega}|}(1)-
\frac{1}{2B^{|\omega|}}\psi\Bigl(\frac{1}{B^{|\omega|}}\Bigr)-
\frac{\gamma}{2B^{|\omega|}}-\frac{1}{2} & \quad \text{if}\quad
v_B(\omega)= 0.
 \end{cases}
 \end{split}
\end{equation*}
\end{corollary}
{\bf Proof.} The required statement easily follows from Theorem
\ref{t5}, Corollary \ref{c4} and the equality
$$
\int_0^1\frac{x^{b-1}(1-x)}{1-x^a}\,dx=\sum_{k=0}^{\infty}
\left(\frac{1}{ak+b}-\frac{1}{ak+b+1}\right)=\frac{1}{a}
\left(\psi\Bigl(\frac{b+1}{a}\Bigr)-\psi\Bigl(\frac{b}{a}\Bigr)\right).
\qed
$$
From Theorem \ref{t3}, (\ref{eq25}) and (\ref{re}) with $a=1,$
$l=0$ we have
\begin{corollary} \label{c5}
Let $b, B>1$ be positive integers, $z\in {\mathbb C},$ $|z|\le 1.$
Then
$$
\gamma_{1,b}(z)=\sum_{k=1}^{\infty}a_kQ(k,B)=\sum_{k=1}^{\infty}a_{\lfloor\frac{k}{B}\rfloor}
\frac{\varepsilon(k)}{k},
$$
where $a_0=a_1=\ldots=a_{b-1}=0,$
$a_k=a_{\lfloor\frac{k}{B}\rfloor}+z^{k-b},$ $k\ge b.$
\end{corollary}

Similarly, from Theorem \ref{t5} we have
\begin{corollary} \label{c65}
Let $b,B>1$ be positive integers, $z\in {\mathbb C},$ $|z|\le 1.$
Then
$$
\gamma_{1,b}(z)=\frac{1}{2}\sum_{k=0}^{\infty}\frac{z^k}{(k+b)(k+b+1)}
+\sum_{k=1}^{\infty}a_k\frac{P_B(k)}{Bk(Bk+1)\cdots (Bk+B)},
$$
where $a_0=a_1=\ldots=a_{b-1}=0,$
$a_k=a_{\lfloor\frac{k}{B}\rfloor}+z^{k-b},$ $k\ge b.$
\end{corollary}

{\bf Example 2.} If in Corollary \ref{c5} we take $z=1,$ then we
get that  $a_k$ is equal to the   $B$-ary length of
$\lfloor\frac{k}{b}\rfloor,$ i. e.,
$$
a_k=\sum_{\alpha=0}^{B-1}N_{\alpha,B}\left(\Bigl\lfloor\frac{k}{b}\Bigr\rfloor\right)
=L_B\left(\Bigl\lfloor\frac{k}{b}\Bigr\rfloor\right).
$$
On the other hand,
$$
\gamma_{1,b}(1)=\log b-\psi(b)=\log
b-\sum_{k=1}^{b-1}\frac{1}{k}+\gamma
$$
and hence we get
\begin{equation}
\log
b-\psi(b)=\sum_{k=1}^{\infty}L_B\left(\Bigl\lfloor\frac{k}{b}\Bigr\rfloor
\right)Q(k,B). \label{eq29}
\end{equation}
If $b=1,$ formula (\ref{eq29}) gives (\ref{eq6}). If $b>1,$ then
from (\ref{eq29}) and (\ref{eq6}) we get
\begin{equation}
\log
b=\sum_{k=1}^{b-1}\frac{1}{k}+\sum_{k=1}^{\infty}\left(L_B\left(\Bigl\lfloor
\frac{k}{b}\Bigr\rfloor\right)-L_B(k)\right)Q(k,B), \label{eq30}
\end{equation}
which is equivalent to \cite[Theorem 2.8]{bebo}. Similarly, from
Corollary \ref{c65} we obtain (\ref{eu}) and
\begin{equation}
\log
b=\sum_{k=1}^{b-1}\frac{1}{k}-\frac{b-1}{2b}+\sum_{k=1}^{\infty}
\frac{\left(L_B(\lfloor
\frac{k}{b}\rfloor)-L_B(k)\right)P_B(k)}%
{Bk(Bk+1)\cdots (Bk+B)}. \label{eq30.5}
\end{equation}

{\bf Example 3.} Using the fact that for any integer $B>1$
$$
L_B\left(\Bigl\lfloor\frac{k}{B}\Bigr\rfloor\right)-L_B(k)=-1,
$$
from (\ref{eq28}), (\ref{eq6}) and (\ref{eq30}) we get the
following rational series for $\log\Gamma(1/B):$
$$
\log\Gamma\left(\frac{1}{B}\right)=\sum_{k=1}^{B-1}\frac{1}{k}
+\sum_{k=1}^{\infty}\Bigl(N_{0,B}(k)-\frac{1}{B}L_B(k)-1\Bigr)Q(k,B).
$$

{\bf Example 4.} Substituting $b=1,$ $z=-1$ in Corollary \ref{c5}
we get the generalized Vacca series for $\log\frac{4}{\pi}.$
\begin{corollary} \label{c6}
 Let $B\in {\mathbb N},$ $B>1.$  Then
$$
\log\frac{4}{\pi}=\sum_{k=1}^{\infty}a_kQ(k,B)=\sum_{k=1}^{\infty}a_{\lfloor\frac{k}{B}\rfloor}
\frac{\varepsilon(k)}{k},
$$
where
\begin{equation}
a_0=0, \qquad a_k=a_{\lfloor\frac{k}{B}\rfloor}+(-1)^{k-1}, \quad
k\ge 1. \label{eq31}
\end{equation}
In particular, if $B$ is even, then
\begin{equation}
\log\frac{4}{\pi}=\sum_{k=1}^{\infty}(N_{odd,B}(k)-N_{even,B}(k))Q(k,B)
=\sum_{k=1}^{\infty}
\frac{\left(N_{odd,B}(\lfloor\frac{k}{B}\rfloor)-
N_{even,B}(\lfloor\frac{k}{B}\rfloor)\right)}{k} \varepsilon(k),
\label{eq32}
\end{equation}
where $N_{odd,B}(k)$ (respectively $N_{even,B}(k)$) is the number
of occurrences  of the odd (respectively even) digits in the
$B$-ary expansion of $k.$
\end{corollary}
{\bf Proof.} To prove (\ref{eq32}), we notice that if $B$ is even,
then the sequence $\widetilde{a}_k:=N_{odd,B}(k)-N_{even,B}(k)$
satisfies  recurrence (\ref{eq31}).  \qed

Substituting $b=1, z=-1$ in Corollary \ref{c65} with the help of
(\ref{eq30.5}) we get the generalized Addison series for
$\log\frac{4}{\pi}.$
\begin{corollary} \label{c69}
Let $B>1$ be a positive integer. Then
$$
\log\frac{4}{\pi}=\frac{1}{4}+\sum_{k=1}^{\infty}\frac{\left(L_B(\lfloor
\frac{k}{2}\rfloor)-L_B(k)+a_k\right)P_B(k)}{Bk(Bk+1)\cdots
(Bk+B)},
$$
where the sequence $a_k$ is defined in Corollary {\rm\ref{c6}}. In
particular, if $B$ is even, then
$$
\log\frac{4}{\pi}=\frac{1}{4}+\sum_{k=1}^{\infty}\frac{\left(L_B(\lfloor
\frac{k}{2}\rfloor)-2N_{even,B}(k)\right)P_B(k)}{Bk(Bk+1)\cdots
(Bk+B)}.
$$
\end{corollary}
{\bf Example 5.} For $t>1,$ the generalized Somos constant
$\sigma_t$ is defined by
$$
\sigma_t=\sqrt[t]{1\sqrt[t]{2\sqrt[t]{3\ldots}}}=
1^{1/t}2^{1/t^2}3^{1/t^3}\cdots=\prod_{n=1}^{\infty}n^{1/t^n}.
$$
In view of the relation \cite[Th.8]{sh}
\begin{equation}
\gamma_{1,1}\left(\frac{1}{t}\right)=t\log\frac{t}%
{(t-1)\sigma_t^{t-1}}, \label{somos}
\end{equation}
by Corollary  \ref{c5} and formula (\ref{eq30}) we get
\begin{corollary} \label{c7}
Let $B\in {\mathbb N},$ $B>1,$ $t\in {\mathbb R},$ $t>1.$ Then
$$
\log\sigma_t=\frac{1}{(t-1)^2}+\frac{1}{t-1}\sum_{k=1}^{\infty}
\left(L_B\Bigl(\Bigl\lfloor\frac{k}{t}\Bigr\rfloor\Bigr)-L_B\Bigl(\Bigl\lfloor
\frac{k}{t-1}\Bigr\rfloor\Bigr)-\frac{a_k}{t}\right)Q(k,B),
$$
where $a_0=0,$ $a_k=a_{\lfloor\frac{k}{B}\rfloor}+t^{1-k},$ $k\ge
1.$
\end{corollary}
In particular, setting $B=t=2$ we get the following rational
series for Somos's quadratic recurrence constant:
$$
\log\sigma_2=1-\frac{1}{2}\sum_{k=1}^{\infty}\frac{a_k}{2k(2k+1)},
$$
where $a_1=3,$
$a_k=a_{\lfloor\frac{k}{2}\rfloor}+\frac{1}{2^{k-1}},$ $k\ge 2.$

From (\ref{somos}), (\ref{eq30.5}) and Theorem \ref{t5} we find
\begin{corollary} \label{c8}
Let $B\in {\mathbb N},$ $B>1,$ $t\in {\mathbb R},$ $t>1.$ Then
\begin{equation*}
\begin{split}
&\log\sigma_t=\frac{3t-1}{4t(t-1)^2} \\
&+\frac{t+1}{2(t-1)}\sum_{k=1}^{\infty}
\left(L_B\Bigl(\Bigl\lfloor\frac{k}{t}\Bigr\rfloor\Bigr)-L_B\Bigl(\Bigl\lfloor
\frac{k}{t-1}\Bigr\rfloor\Bigr)-\frac{2a_k}{t(t+1)}\right)
\frac{P_B(k)}{Bk(Bk+1)\cdots(Bk+B)},
\end{split}
\end{equation*}
 where the sequence $a_k$ is
defined in Corollary {\rm \ref{c7}}.
\end{corollary}
In particular, if $B=t=2$ we get
$$
\log\sigma_2=\frac{5}{8}-\frac{1}{2}\sum_{k=1}^{\infty}
\frac{a_k}{2k(2k+1)(2k+2)},
$$
where $a_1=4,$
$a_k=a_{\lfloor\frac{k}{2}\rfloor}+\frac{1}{2^{k-1}},$ $k\ge 2.$

{\bf Example 6.} The Glaisher-Kinkelin constant is defined by the
limit \cite[p.135]{fi}
$$
A:=\lim_{n\to\infty}\frac{1^22^2\cdots
n^n}{n^{\frac{n^2+n}{2}+\frac{1}{12}}e^{-\frac{n^2}{4}}}=1.28242712\ldots.
$$
Its connection to the generalized-Euler-constant function
$\gamma_{a,b}(z)$ is given by the formula \cite[Cor.4]{sh}
\begin{equation}
\gamma'_{1,1}(-1)=\log\frac{2^{11/6}A^6}{\pi^{3/2}e}.
\label{glaisher}
\end{equation}
By Theorem \ref{t5}, since
$$
\int_0^1\frac{x(1-x)}{(1+x)^2}\,dx=3\log 2-2,
$$
we have
$$
\log A=\frac{4}{9}\log 2-\frac{1}{4}\log\frac{4}{\pi}+
\frac{1}{6}\sum_{k=1}^{\infty}a_{k,1}\frac{P_B(k)}{Bk(Bk+1)\cdots(Bk+B)},
$$
where the sequence $a_{k,1}$ is defined by the generating function
(\ref{eq24}) with $a=b=l=1,$ $z=-1,$ or using (\ref{re}) by the
recursion
$$
a_{0,1}=a_{1,1}=0, \qquad
a_{k,1}=a_{\lfloor\frac{k}{B}\rfloor,1}+(-1)^k(k-1), \quad k\ge 2.
$$
Now by Corollary \ref{c69} and (\ref{eq30.5}) we get
\begin{corollary} \label{c12}
Let $B>1$ be a positive integer. Then
$$
\log A=\frac{13}{48}-\frac{1}{36}\sum_{k=1}^{\infty}
\left(7L_B(k)-7L_B\Bigl(\Bigl\lfloor\frac{k}{2}\Bigr\rfloor\Bigr)
+b_k\right)\frac{P_B(k)}{Bk(Bk+1)\cdots(Bk+B)},
$$
where $b_0=0,$
$b_k=b_{\lfloor\frac{k}{B}\rfloor}+(-1)^{k-1}(6k+3),$ $k\ge 1.$
\end{corollary}
In particular, if $B=2$ we get
$$
\log A=\frac{13}{48}-\frac{1}{36}\sum_{k=1}^{\infty}
\frac{c_k}{2k(2k+1)(2k+2)},
$$
where $c_1=16,$
$c_k=c_{\lfloor\frac{k}{2}\rfloor}+(-1)^{k-1}(6k+3),$ $k\ge 2.$

Using the formula expressing $\frac{\zeta'(2)}{\pi^2}$ in terms of
Glaisher-Kinkelin's constant \cite[p.135]{fi}
$$
\log A=-\frac{\zeta'(2)}{\pi^2}+\frac{\log 2\pi+\gamma}{12}
$$
by Corollaries \ref{c65}, \ref{c69} and \ref{c12}, we get
\begin{corollary}
Let $B>1$ be a positive integer. Then
$$
\frac{\zeta'(2)}{\pi^2}=-\frac{1}{16}+\frac{1}{36}\sum_{k=1}^{\infty}
\left(4L_B(k)-L_B\Bigl(\Bigl\lfloor\frac{k}{2}\Bigr\rfloor\Bigr)+c_k\right)
\frac{P_B(k)}{Bk(Bk+1)\cdots(Bk+B)},
$$
where $c_0=0,$ $c_k=c_{\lfloor\frac{k}{B}\rfloor}+(-1)^{k-1}6k,$
$k\ge 1.$
\end{corollary}

{\bf Example 7.}  First we evaluate $\gamma_{2,1}^{(l)}(-1)$ for
$l=0, 1.$ From Corollaries  \ref{c1}, \ref{c1.1} and
\cite[Ex.3.12, 3.13]{sh} we have
$$
\gamma_{2,1}(-1)=\int_0^1\int_0^1\frac{(x-1)\,dxdy}{(1+x^2y^2)\log
xy}=\frac{\pi}{4}-2\log\Gamma\Bigl(\frac{1}{4}\Bigr)+\log\sqrt{2\pi^3}
$$
and
\begin{equation*}
\begin{split}
\gamma'_{2,1}(-1)=&-\frac{1}{4}\Phi(-1,1,3/2)+\frac{1}{2}
\Phi(-1,0,3/2) +\frac{1}{2}\frac{\partial\Phi}{\partial
s}(-1,0,3/2) \\
&-\frac{\partial\Phi}{\partial s}(-1,-1,3/2)-
\frac{\partial\Phi}{\partial
s}(-1,0,2)+\frac{\partial\Phi}{\partial s}(-1,-1,2).
\end{split}
\end{equation*}
The last expression
can be evaluated explicitly (see \cite[Section 2]{sh}) and we get
$$
\gamma'_{2,1}(-1)=\frac{G}{\pi}+\frac{\pi}{8}-\log\Gamma\Bigl(
\frac{1}{4}\Bigr)-3\log A+\log\pi+\frac{1}{3}\log 2,
$$
or
\begin{equation}
\frac{G}{\pi}=\gamma'_{2,1}(-1)-\frac{1}{2}\gamma_{2,1}(-1)
+\frac{1}{4}\log\frac{4}{\pi}+3\log A-\frac{7}{12}\log 2.
\label{eq101}
\end{equation}
On the other hand, by Theorem \ref{t5} and (\ref{re}) we have
\begin{equation}
\gamma_{2,1}(-1)=\frac{\pi}{8}-\frac{1}{4}\log
2+\sum_{k=1}^{\infty}a_{k,0}\frac{P_B(k)}{Bk(Bk+1)\cdots (Bk+B)},
\label{eq102}
\end{equation}
where $a_{0,0}=0,$ $a_{2k,0}=a_{\lfloor\frac{2k}{B}\rfloor,0},$
$k\ge 1,$ $a_{2k+1,0}=a_{\lfloor\frac{2k+1}{B}\rfloor,0}+(-1)^k,$
$k\ge 0,$ and
\begin{equation}
\gamma'_{2,1}(-1)=\frac{\pi}{16}-\frac{1}{4}\log
2+\sum_{k=1}^{\infty}a_{k,1}\frac{P_B(k)}{Bk(Bk+1)\cdots (Bk+B)},
\label{eq103}
\end{equation}
where $a_{0,1}=0,$ $a_{2k,1}=a_{\lfloor\frac{2k}{B}\rfloor,1},$
$k\ge 1,$
$a_{2k+1,1}=a_{\lfloor\frac{2k+1}{B}\rfloor,1}+(-1)^{k-1}k,$ $k\ge
0.$ Now from (\ref{eq101}) -- (\ref{eq103}), (\ref{eq30.5}) and
Corollary \ref{c69} we get the following expansion for $G/\pi.$
\begin{corollary} Let $B>1$ be a positive integer. Then
$$
\frac{G}{\pi}=\frac{11}{32}+\sum_{k=1}^{\infty}\left(
\frac{1}{8}L_B\Bigl(\Bigl\lfloor\frac{k}{2}\Bigr\rfloor\Bigr)
-\frac{1}{8}L_B(k)+c_k\right)\frac{P_B(k)}{Bk(Bk+1)\cdots (Bk+B)},
$$
where $c_0=0,$ $c_{2k}=c_{\lfloor\frac{2k}{B}\rfloor}+k,$ $k\ge
1,$ $c_{2k+1}=c_{\lfloor\frac{2k+1}{B}\rfloor}+
\frac{(-1)^{k-1}-1}{2}(2k+1),$ $k\ge 0.$
\end{corollary}
In particular, if $B=2$ we get
$$
\frac{G}{\pi}=\frac{11}{32}+\sum_{k=1}^{\infty}\frac{c_k}{2k(2k+1)(2k+2)},
$$
where $c_1=-\frac{9}{8},$ $c_{2k}=c_k+k,$ $c_{2k+1}=c_k+
\frac{(-1)^{k-1}-1}{2}(2k+1),$ $k\ge 1.$

\vspace{1cm}

\end{document}